\documentclass{article}
\usepackage{graphicx}
\usepackage{amsmath}
\usepackage{amsfonts}
\usepackage{amssymb}
\newcommand{\RM}{\mathbb{R}}
\newcommand{\ZM}{\mathbb{Z}}

\newcommand{\real}{{\rm Re}}
\newcommand{\imag}{{\rm Im}}

\newcommand{\bx}{\mathbf{x}}
\newcommand{\by}{\mathbf{y}}

\newcommand{\bv}{\mathbf{v}}
\newcommand{\bn}{\mathbf{n}}

\newcommand{\curl}{\operatorname{curl}}
\newcommand{\sign}{\operatorname{sign}}
\newcommand{\ds}{\operatorname{ds}}
\newcommand{\dl}{\operatorname{dl}}

\numberwithin{equation}{section}

\begin{document}

\title{ENSTROPHY DYNAMICS FOR FLOW PAST A SOLID BODY  WITH NO-SLIP BOUNDARY CONDITION}

\author{Aleksei Gorshkov}


\maketitle

\begin{abstract}
In the paper we study the impact of the boundary vorticity distribution on the dynamics of enstrophy for flows around streamlined body. A new energy identity is derived in the article, which includes the boundary values of the vortex function. For the Stokes system the dissipativity of enstrophy is proved. For the Navier-Stokes system a new equation of the enstrophy dynamics is obtained.

\smallskip

{\it keywords:} Enstrophy, Vorticity equation, dissipation, energy identity.
\end{abstract}

\bigskip
The article examines the dynamics of {\it enstrophy} in flows around a solid body with the no-slip condition at its boundary. The enstrophy is defined as the square of the $L_2$-norm of the vortex function $w$:
$$\mathcal E(t) = \int_\Omega w^2(t,\bx) d\bx.$$  

Under periodic boundary conditions or in the case of empty boundary, the dynamics of enstrophy is expressed via the {\it palinstrophy} by the equation:
\begin{align*}
\mathcal E(t) = - 2 \nu \int_\Omega |\nabla w(\bx)|^2d\bx,
\end{align*}
ensuring its decrease at $t \to \infty$ (here $\nu$ is a viscosity). The dissipation of enstrophy along with the energy cascade for periodic solutions was studied in \cite{BaiesiMaes}.  The behavior of dissipation at $\nu \to 0$ in the periodic case was investigated in \cite{MPY}. However, in the problem of external flow around a body, the equation above is lost due to the presence of boundary conditions.

The interaction of a Newtonian fluid with a solid in terms of vorticity has been the subject of many studies, and despite a number of fundamental results and general theories in this field, this problem still raises many open questions (see \cite{Wu} and the literature review cited there). And one of the problems in this context is the distribution of vorticity near the body boundary, given that the source of vorticity is both non-viscosity components (boundary velocity distribution, pressure gradient) and viscosity-dependent tangential stresses at the boundary. The type of equations of motion determines the nature of the generated vorticity and its main source. In Euler's equations, these are "inviscid" sources, while for the Stokes and Navier-Stokes systems, dissipative terms related to friction play a significant (but not exclusive) role. In this paper, the equation of the dynamics of enstrophy will be presented with the division of terms into boundary and distributed, sign-defined and alternating, viscous and inviscid.

The evolution of enstrophy and its influence on flow characteristics has been studied by many authors. The relationship of enstrophy with the strain tensor was investigated by Weiss in \cite{Weiss}. The works of Lighthill\cite{Lighthill}, Wu\cite{Wu2}, Chen, Li, Wang\cite{Chen} investigated the evolution of the normal derivative of enstrophy at the body boundary, and its relationship with viscous friction and pressure in the boundary layer.

\section{Energy identities for a vector field}
The dynamics of the vortex function $w$ under the assumption of the absence of external forces is described by the Helmholtz vortex equation. The derivative of enstrophy satisfies the equation:
\begin{align}\label{enstrdyn}
\frac 12 \frac d{dt} \int_\Omega w^2(t,\bx)d\bx + \nu \int_\Omega |\nabla w|^2 d\bx - \nu \oint \limits_{\partial \Omega} \frac {\partial w}{\partial \bn} w (t,\bx) \dl =0.
\end{align}

Here $\nu$ is the coefficient of kinematic viscosity. The positive term $\int_\Omega |\nabla w|^2 d\bx$ characterizes the viscous dissipation of enstrophy. It defines the enstrophy cascade in transition from large to smaller scales. When the no-slip condition holds, the term
\begin{align}\label{boundenstr}
\oint \limits_{\partial \Omega} \frac {\partial w}{\partial \bn} w (t,\bx)\dl
\end{align}
obtains a different representation based on the boundary vorticity distribution. In fact, this value characterizes the additional contribution of the no-slip condition to the process of generating vortex motion and 
in work done to overcome viscous forces around the streamlined body. For the Stokes system, the dynamics of enstrophy will be determined by two quantities opposite in sign - dissipation and boundary vorticity caused by the no-slip condition. And the contribution of the no-slip condition to the dynamics of enstrophy will be quantified in terms of the Sobolev norm with a fractional exponent of 1/2 of vorticity  $w(t,\bx)$ distributed along the boundary. For a nonlinear Navier-Stokes system with a zero velocity on the streamlined body, the value in (\ref{boundenstr}) will no longer be sign-defined. For well-streamlined bodies, the boundary value is $w(t,\bx)$ will continue to make a positive contribution to the enstrophy dynamics. But with increasing vorticity, its dynamics will become more uncertain.

The importance of studying the dynamics of enstrophy is also caused by its connection with the problem of  the regularity and uniqueness of the solution solution to the Navier-Stokes system. The boundedness of the enstrophy as a function of time allows us to prove the smoothness of the solution of the Navier-Stokes system. Uniform estimates of $\mathcal E(t)$ ensure both the existence and uniqueness of the solution \cite{FoiasTemam}.

Turbulent flows are characterized by a large fluctuations of vorticity. For incompressible flows with finite kinetic energy and zero boundary condition, the equality holds $\|\nabla \bv\|_{L_2} = \|w\|_{L_2}$ (Doering, Gibbon \cite{DoeringGibbon}), and so the enstrophy $\mathcal E$ determines the rate 
the rate of kinetic energy dissipation:
\begin{align}\label{energydissip}
\frac 12 \frac d{dt}\int_\Omega |\bv(\bx)|^2 d\bx = - \nu \int_\Omega |w(\bx)|^2d\bx.
\end{align}
Enstrophy is also used in turbulent models for a number of near-wall and other types of flows. Thus, sometimes in $k-\epsilon$ turbulence models, instead of the equation for energy dissipation $\epsilon$, the corresponding equation for enstrophy is used. For some flows it describes the process more adequately \cite{TL}. 

The rate of internal energy of a viscous medium is conditioned by the work of dissipative forces done in deforming the element of fluid. It is determined by the deviatoric stress and the non-isotropic part of the strain tensor\cite{Batchelor}. For incompressible flows with no-slip condition, the deviatoric stress can be found from formula:
\begin{align}\label{straincurl}
{\frac {\nu}{2}}\int_\Omega \left({\frac {\partial u_{i}}{\partial x_{j}}}+{\frac {\partial u_{j}}{\partial x_{i}}}\right)^2 d\bx = \nu \mathcal E.
\end{align}
Moreover, the relationship between dissipation and enstrophy is exclusively non-local and involves $L_2$-
norm of these values calculated by the entire region. In the work \cite{ZhuAntonia}, the correlation of local values of the characteristics of enstrophy and energy was investigated. If the no-slip condition is replaced by the slip condition (the normal velocity component at the boundary is zero), then the formula (\ref{straincurl}) is supplemented by a number of boundary terms. And the sign of the difference between the strain tensor and the enstrophy begins to determine the type of the flow (see \cite{Weiss}). The relationship of vorticity and viscous dissipation in the flows around rotating bodies in the incompressible medium was studied in \cite{Koh}. The dissipation of the vorticity itself, including experimental studies, was studied by Karman\cite{Karman} and Taylor\cite{Taylor}.


\vskip 10pt
The article explores the two-dimensional external problem of fluid dynamics that describes a flow which interacts with body 
by no-slip condition. A new equation of vorticity energy dynamics $\mathcal E(t)$ will be presented, in which the value  (\ref{boundenstr}) obtained from the Green formula will decompose into sign-defined terms. Thus, the rate of enstrophy  will directly depend on the distribution of vorticity along the boundary of the streamlined body.

Energy estimates make sense only if the energy of the flow is finite, which is often not the case in unbounded areas.
Flows around the solids are the ones with infinite energy due to the presence of a nonzero boundary condition at infinity, namely the velocity of the incoming flow $\bv_\infty$.
However, even after subtracting $\bv_\infty$ from the flow, it does not always have a finite energy. So, if, for example, the circulation of the current at infinity is not zero, then $$E=\int |\bv-\bv_\infty|^2 d\bx = \infty.$$

For example, the delta function of a discrete vortex on a plane generates a circulation flow
\begin{align}\label{pointvortex}
\bv(\bx) =  \bv_\infty + \frac {\bx^\perp} {|\bx|^2},
\end{align}
which is not quadratically integrable even when subtracting $\bv_\infty$. If the flow around the streamlined body
with no-slip condition at its boundary has a non-zero mean curl value(or divergence in the case of compressible fields), then such a flow will also have an infinite norm $\|\bv-\bv_\infty\|_{L_2}$. This follows from the Stokes and Gauss-Ostrogradsky theorems, since nonzero mean values of the curl generate nonzero circulation flows at infinity of the form (\ref{pointvortex}). In the paper of the author \cite{G2}, it is proved that if the field $\bv$ is solenoidal, satisfies the no-slip condition, $\curl\bv$ is quadratically integrable, and
\begin{align} \label{intro:zerocirculation}
&\lim_{R\to\infty}\oint_{\vert\bx\vert=R} \bv \cdot d\mathbf{l} = 0,
\end{align}  
then $\| \bv-\bv_\infty \|^2_{L_2(\Omega)} < \infty$.

The article considers flows with a finite $L_2$-norm of the vector field $\bv-\bv_\infty$, for which  (\ref{intro:zerocirculation}) holds. For them, instead of the identity (\ref{energydissip}), the following equality is valid
\begin{align*}
\frac 12 \frac d{dt}\int_\Omega |\bv(\bx) - \bv_\infty|^2 d\bx = - \nu \int_\Omega |w(\bx)|^2d\bx.
\end{align*}
In turn, the energy identity for enstrophy obtained in this article will already describe the dynamics of the rate of energy change, namely $\frac{d^2}{dt^2}\int_\Omega |\bv(\bx)|^2 d\bx$.

\section{Dissipation of the enstrophy of the Stokes system}
\subsection{Energy identity in the exterior of the disk} 
The curl $$w(t,\bx)={\rm curl}~\bv(t,\bx)=\partial_{\bx_1}v_2 -  \partial_{\bx_2}v_1$$ 
for 2D vector fields $\bv(t,\bx)=(v_1(t,\bx),v_2(t,\bx))$, described by the Stokes system, satisfies the heat equation:
\begin{align}
\frac{\partial w(t,\bx)}{\partial t}	 - \nu \Delta w = 0.  \label{heateq} 
\end{align}

In exterior of the disk $B_{r_0}=\{\bx \in \RM^2,~|\bx| > r_0 \},~r_0>0$ with a fixed horizontal flow at infinity $\bv_\infty = (v_{\infty},0) \in \RM^2$ the no-sip condition
\begin{align}\label{noslip}
\bv(\bx)=0,~|\bx| = r_0
\end{align}
is the orthogonality condition for vorticity(see \cite{G2}):
\begin{align}\label{ortcond_}
\frac 1{2\pi} \int_{\Omega} \frac {w(t,\bx)}{z^k} d\bx = \begin{cases} 0,~k\in \mathbb{N}\cup \{0\},~k\neq 1, \\ i v_\infty,~k=1.
\end{cases}
\end{align}

For the Stokes flow, this condition becomes a boundary ones, and is given by the relations:
\begin{equation}\label{robin_bound}
r_0\frac{\partial w_k(t,r)}{\partial r}\Big|_{r=r_0} + |k| w_k(t,r_0) = 0,~k \in \ZM.
\end{equation}
Here $w_k(t,r)$ are the Fourier coefficients of the function $w$ decomposed by $\{e^{ik\varphi}\}$.

From the Parseval equality, in virtue of (\ref{robin_bound}), and $ \frac {\partial}{\partial \bn} = -\frac {\partial }{\partial r}$, we have
\begin{align*}
\oint \limits_{|\bx| = r_0} \frac {\partial w}{\partial \bn} w(t,\bx) \dl = -r_0 \sum_{k=-\infty}^\infty
\frac{\partial w_k(t,r)}{\partial r}\Big|_{r=r_0} w_k(t,r_0) = 
\sum_{k=-\infty}^\infty
|k| w^2_k(t,r).
\end{align*}

On the circle $S_{r_0}=\{\bx \in \RM^2,~|\bx| = r_0 \}$ the Sobolev semi-norm with a fractional exponent $\frac 12$:
$$
\dot H^{\frac 12}(S_{r_0}) = \left \{f(\cdot) \in L_2(S_{r_0}),~ \sum_{k=-\infty}^\infty |k|f_k^2 <\infty ~\right \},
$$
where $f_k$ are the Fourier coefficients of the function $f$.

As a result, we obtain a new identity describing the dynamics of the vorticity energy:
\begin{align*}
\frac 12 \frac d{dt} \|w(t,\cdot)\|^2_{L_2(B_{r_0})} + \nu \|\nabla w\|^2_{L_2(B_{r_0})} - \nu \|w\|^2_{\dot H^{\frac 12}(S_{r_0})} =0
\end{align*}

It is well known that $H^{\frac 12}$-norm is controled by Sobolev norm $H^1$ (see \cite{Lions}). It will be proved below that for solutions of the Stokes system with the no-slip condition the following is true: $$\|w\|^2_{\dot H^{\frac 12}(\partial \Omega)}  \leq \|\nabla w\|^2_{L_2(\Omega)}.$$

\subsection{Energy identity in the exterior of a simply connected area.}
Now we consider the case of an external domain satisfying certain conditions. Let the domain $\Omega=\RM^2\setminus\overline G$ be the exterior of a bounded simply connected domain $G$ the boundary of which is a smooth Jordan curve. We will also assume that there is a Riemann mapping
$$\Phi:\Omega \to B_{r_0}$$ 
such that
\begin{align*}
\Phi(p)=p+O\left (\frac 1p \right ),
\end{align*}
and the inverse mapping $\Phi^{-1}(z) : B_{r_0} \to \Omega$ satisfies:
\begin{align*}
&\Phi^{-1}(z)=z+O\left (\frac 1z \right ),\\
&\left(\Phi^{-1}\right )'(z)=1+O\left (\frac 1{z^2} \right ).
\end{align*}

The Riemann mapping will be applied both to the complex variable $z$ and to the corresponding real vector $\bx= (\real z,\imag z)$ retaining the same designation $\Phi$ for both cases. 

For a flow satisfying the boundary condition at infinity
$$
\bv(\bx)\to \bv_\infty,~|\bx|\to \infty
$$
with given velocity $\bv_\infty=(v_\infty, 0)$,
the no-slip condition (\ref{noslip}) is equivalent to the orthogonality condition(see \cite{G2}):
\begin{align}\label{ortcond}
\frac 1{2\pi} \int_{\Omega} \frac {w(t,\bx)}{\Phi(z)^k} d\bx = \begin{cases} 0,~k\in \mathbb{N}\cup \{0\},~k\neq 1, \\ i v_\infty,~k=1.
\end{cases}
\end{align}

Change the variables $p=y_1+iy_2 \in \Omega$ by
$z=x_1+ix_2 \in B_{r_0}$ in equation (\ref{heateq}) where $z=\Phi(p)$. The Laplace operator $\Delta$ transferes to $\vert \Phi'(p) \vert^2 \Delta$ and the heat equation takes the form
$$
\frac{\partial w(t,\bx)}{\partial t} - \nu \vert \Phi'(p) \vert^2 \Delta w=0,
$$
or
\begin{align}\label{heateq2}
\vert (\Phi^{-1})'(z) \vert^2 \frac{\partial w(t,\bx)}{\partial t} - \nu \Delta w=0.
\end{align}
Here for the curl $w\left (t, \Phi^{-1}(\bx) \right)$, defined in $B_{r_0}$ we have retained the same designation $w$.

Fix integer $k \geq 0$ and divide the equation (\ref{heateq2}) by $z^k=(x_1+ix_2)^k$. Integrate it over $B_{r_0}$. From (\ref{ortcond}) follows
\begin{align}\label{ortcond2}
\frac d{dt}\int_{B_{r_0}} \frac {\vert (\Phi^{-1})'(z) \vert^2}{z^k} w(t,\bx) d\bx = 0.
\end{align}

Then
$$
\int_{B_{r_0}}  \frac {\Delta w}{z^{k}}d\bx=0.
$$

On the other hand
\begin{align*}
\int_{B_{r_0}} \frac {\Delta w}{z^{k}}d\bx = \int_{r_0}^\infty \int_0^{2\pi} \frac {\Delta w}{s^ke^{ik\varphi}}sdsd\varphi = 2\pi \int_{r_0}^\infty s^{-k+1}\Delta_k w_k(s) ds \\= -2\pi r_0^{-k}\left(r_0 \frac{\partial w_k(t,r)}{\partial r}\Big \vert_{r=r_0} + k w_k(t,r_0) \right )=0.
\end{align*}

From $w_{-k}(t,r) = \overline{w_k(t,r)}$, we obtain for all integer $k$ the same boundary condition (\ref{robin_bound}), valid in an arbitrary external domain satisfying the constraints imposed above.

Multiply  (\ref{heateq2}) by $w$ and integrate it over $B_{r_0}$:
\begin{align*}
\frac 12 \frac d{dt} \int_{B_{r_0}} w^2(t,\bx) \vert (\Phi^{-1})'(z) \vert^2 d\bx + \nu \int_{B_{r_0}} |\nabla w|^2 d\bx - \nu \oint \limits_{|\bx|=r_0} \frac {\partial w}{\partial \bn} w (t,\bx) \dl =0.
\end{align*}

From the boundary condition (\ref{robin_bound}) we will have
\begin{align*}
\frac 12 \frac d{dt} \int_{B_{r_0}} w^2(t,\bx) \vert (\Phi^{-1})'(z) \vert^2 d\bx + \nu \int_{B_{r_0}} |\nabla w|^2 d\bx - \nu \sum_{k=-\infty}^\infty
|k| w^2_k(t,r) =0.
\end{align*}

Define a seminorm of the Sobolev space with a fractional exponent:
$$\|f\|_{\dot H^{\frac 12}(\partial \Omega)} = \|f \circ \Phi^{-1}\|_{\dot H^{\frac 12}(S_{r_0})}.$$

Then we finally get the energy identity:
\begin{align}\label{energy1}
\frac 12 \frac d{dt} \int_\Omega w^2(t,\bx)d\bx + \nu \int_\Omega |\nabla w|^2 d\bx - \nu \|w(t,\cdot)\|^2_{\dot H^{\frac 12}(\partial \Omega)} =0.
\end{align}

\vskip 10pt
From this identity follows that the dissipation of vorticity decreases due to viscous forces distributed over the external region, and at the same time increases due to vorticity distributed at the boundary of the body. Viscous forces and boundary vorticity are of opposite sign in contribution to the dynamics of the enstrophy. $H^{\frac 12}$-norm of vorticity characterizes an increase in dissipation caused by the slowdown of the flow around the streamlined body. Nevertheless, we will prove that the distributed viscous forces prevail over the boundary ones, ensuring that the curl stabilizes to zero.

\subsection{Enstrophy dissipation}
Let's prove the strict negativity of the quadratic form generated by the operator $\Delta$ with boundary condition (\ref{noslip}):
\begin{align*}
(\Delta w,w) =  - \|\nabla w\|^2_{L_2(B_{r_0})} + \|w\|^2_{\dot H^{\frac 12}(S_{r_0})} <0.
\end{align*}
In author's paper \cite{G1} was studied the spectrum of $\Delta$ with a boundary condition (\ref{robin_bound}). The spectrum consists of a continuous part $(-\infty,0)$ and an eigen value $\{0\}$ from the kernel of $\Delta$. From the orthogonality conditions (\ref{ortcond}) follows the orthogonality of the vortex function $w$ to zero eigen value, and the spectral representation of the solution involves only the negative values of the continuous spectrum $\Delta$. The last fact causes strict decrease of enstrophy.

We will find a representation of the quadratic form $-(\Delta w,w)$ as a sum of squares of some differential operators. To do this, we use the decomposition of the Laplace operator into a Fourier series with the operators
$$
\Delta_k f_k(r) = \frac 1r \frac {d}{d r}\left(r \frac {d f_k}{d r}\right) - \frac{k^2}{r^2} f_k(r).
$$ 

Consider the form generated by the operator $\Delta_k$:
$$
(\Delta_k f, f) = \int_{r_0}^\infty \left ( \frac 1r \frac {d}{d r}\left(r \frac {d f}{d r}\right) - \frac{k^2}{r^2} f(r) \right ) f(r) rdr.
$$

After integration in parts, we obtain the following expression:
\begin{align*}
(\Delta_k f, f) = -\int_{r_0}^\infty \left ( \left ( \frac {d f}{d r} \right )^2  + \frac{k^2}{r^2} f^2(r) \right ) rdr - r_0 f(r_0) \frac {d f}{d r} \Bigg |_{r=r_0}.
\end{align*}

From (\ref{robin_bound}) we get the Boltz functional, which has a maximum at zero:
\begin{align*}
(\Delta_k f, f) = -
\int_{r_0}^\infty \left ( \left ( \frac {d f}{d r} \right )^2  + \frac{k^2}{r^2} f^2(r) \right ) rdr + |k| f^2(r_0) = \\-
\int_{r_0}^\infty \left (  \frac {d f}{d r}  + \frac{|k|}{r} f(r) \right ) ^2 rdr < 0.
\end{align*}

Thus, $\Delta_k$ is a dissipative operator, and, according to the Lumer-Phillips theorem, it is a generator of a $C^0$-semigroup.

Let's define the operator
\begin{align}\label{DkOp}
D_k[f] = \left (  \frac {d f}{d r}  + \frac{|k|}{r} f(r) \right ).
\end{align}

Then according to Parseval's equality
$$
(\Delta w,w) = - \sum_{k=-\infty}^\infty \left | D_k[w_k] \right |^2,
$$
and we get another energy identity:
\begin{align}\label{energy2}
\frac 12 \frac d{dt} \|w(t,\cdot) \|^2_{L_2(\Omega)} + \nu \sum_{k=-\infty}^\infty \left | D_k[w_k] \right |^2=0.
\end{align}

\section{Enstrophy Dynamics for the Navier-Stokes system }

For the fluid flows described by nonlinear equations, the dissipative property of enstrophy is lost. Vortex dynamics is described by the Helmholtz equation
\begin{align}
\frac{\partial w(t,\bx)}{\partial t} - \nu \Delta w + (\bv, \nabla w) = 0.  \label{helmeq} 
\end{align}
The energy identity (\ref{enstrdyn}) remains the same, but the mechanism of influence of the boundary vorticity (\ref{boundenstr}) on enstrophy becomes different.

Fix a non-negative integer $k$ and divide the above equation by $z^k = (x_1+ix_2)^k$. From (\ref{ortcond_}) we have
\begin{align*}
\frac d{dt} \int_{\Omega} \frac {w(t,\bx)}{z^k} d\bx = 0. 
\end{align*}

Then
\begin{align*}\int_{\Omega} \frac {\Delta w(t,\bx)}{z^k} d\bx = 
2\pi \int_{r_0}^\infty s^{-|k|+1} \Delta_k w_k(t,s)\ds = \\-2\pi r_0^{-|k|}\left(r_0 \frac{\partial w_k(t,r)}{\partial r}\Big|_{r=r_0} + k w_k(t,r_0) \right ).
\end{align*}


For $k\geq 0$, we obtain a combined no-slip condition consisting of a boundary operator and a condition distributed over the entire domain:
\begin{align*}
2\pi r_0^{-k}\left(r_0 \frac{\partial w_k(t,r)}{\partial r}\Big|_{r=r_0} + k w_k(t,r_0) \right )
+ \int_{\Omega} \frac {(\bv, \nabla w)}{z^k} d\bx =0.
\end{align*}

In virtue of $w_{-k}(t,r) = \overline{w_k(t,r)}$, then, taking the complex conjugation in the previous formula, we obtain a similar condition for negative $k\in\ZM$:
\begin{align*}
2\pi r_0^{-|k|}\left(r_0 \frac{\partial w_k(t,r)}{\partial r}\Big|_{r=r_0} + |k| w_k(t,r_0) \right )
+ \int_{\Omega} \frac {(\bv, \nabla w)}{\overline {z^k} } d\bx =0.
\end{align*}

Fourier coefficients of the Biot-Savart integral operator $$G[w] = \frac 1{2\pi} \int_{B_{r_0}} \frac{(\bx-\by)^\perp}{\vert\bx-\by\vert^2} w(\by) \operatorname{d\by} = \sum_{k=-\infty}^\infty b_{k}(r)e^{ik\varphi},$$ written out for the radial and tangential parts at $|\bx|=r_0$, are given by formulas (see \cite{G2}):
\begin{align*}b_{k}(r)=
\begin{pmatrix}
\sign(k) \frac{ir^{|k|-1}}2 \int_{r_0}^\infty s^{-|k|+1}w_k(s)\ds\\
-\frac{r_0^{|k|-1}}2 \int_{r_0}^\infty s^{-|k|+1}w_k(s)\ds
\end{pmatrix}.
\end{align*}

At the same time, for $k \in \ZM_+$
\begin{align*}
\int_{r_0}^\infty s^{-|k|+1}w_k(s)\ds =\frac 1{2\pi} \int_{B_{r_0}} \frac {w(t,\bx)}{z^k} d\bx,
\end{align*} 
and for $k \in \ZM_-$
\begin{align*}
\int_{r_0}^\infty s^{-|k|+1}w_k(s)\ds =\frac 1{2\pi} \int_{B_{r_0}} \frac {w(t,\bx)}{\overline {z^{|k|}}} d\bx.
\end{align*}

Hence the integrals $$\int_{\Omega} \frac {(\bv, \nabla w)}{z^k} d\bx, \int_{\Omega} \frac {(\bv, \nabla w)}{\overline {z^k}} d\bx$$ can be found through the Fourier coefficients of the tangential or radial part of a bilinear form:
\begin{align*}
\frac 1{2\pi} \int_{B_{r_0}} \frac{(\bx-\by)^\perp}{\vert\bx-\by\vert^2} \left(\bv, \nabla w(\by)\right) \operatorname{d\by} 
\end{align*}
at $|\bx|=r_0$.

The tangential component of the Biot-Savart operator is defined as
$$
G_\tau[w] = \frac 1{2\pi}  \left ( \int_{B_{r_0}} \frac{(\bx-\by)^\perp}{\vert\bx-\by\vert^2} w(\by) \operatorname{d\by}, \vec \tau \right ),
$$
where $\vec \tau$ - tangent vector to the circle $|\bx|=r_0$ turned  counterclockwise.

We obtain the following condition for the curl $w$:
\begin{align*}
r_0 \frac{\partial w_k(t,r)}{\partial r}\Big|_{r=r_0} + |k| w_k(t,r_0)  = 2 r_0 [G_\tau (\bv, \nabla w)]_k \Big |_{r=r_0},
\end{align*}
where $k$ denotes the corresponding Fourier coefficient. Contour integral in (\ref{enstrdyn}) will be found as
\begin{align*}
\oint \limits_{|\bx| = r_0} \frac {\partial w}{\partial \bn} w(t,\bx) \dl = -r_0 \sum_{k=-\infty}^\infty
\frac{\partial w_k(t,r)}{\partial r}\Big|_{r=r_0} w_k(t,r_0) =\\ 
\sum_{k=-\infty}^\infty 
|k| w^2_k(t,r) - 2 \oint \limits_{|\bx| = r_0} w \cdot G_\tau [ (\bv, \nabla w) ]  \dl.
\end{align*}

This integral can be reduced to contour integral of a vector field:
$$
\oint \limits_{|\bx| = r_0} w \cdot G_\tau [ (\bv, \nabla w) ]  \dl = \oint \limits_{|\bx| = r_0} w \cdot G [ (\bv, \nabla w) ]  \operatorname{\mathbf{dl}}.
$$ 

Then the dynamics of enstrophy in the exterior of the disk will be described by the equation
\begin{align*}
\frac 12 \frac d{dt} \|w(t,\cdot)\|^2_{L_2(B_{r_0})} + \nu \|\nabla w\|^2_{L_2(B_{r_0})} - \nu \|w\|^2_{\dot H^{\frac 12}(S_{r_0})} + 2 \oint \limits_{|\bx| = r_0} w \cdot G [ (\bv, \nabla w) ]  \operatorname{\mathbf{dl}} = 0,
\end{align*}
or
\begin{align*}
\frac 12 \frac d{dt} \|w(t,\cdot) \|^2_{L_2(B_{r_0})} + \nu \sum_{k=-\infty}^\infty \left | D_k[w_k] \right |^2 + 2 \oint \limits_{|\bx| = r_0} w \cdot G [ (\bv, \nabla w) ]  \operatorname{\mathbf{dl}} = 0.
\end{align*}

\vskip 10pt
Consider the general case of an external domain $\Omega$, satisfying the conditions imposed in the previous paragraph.
The Riemann mapping $\Phi$ transforms the Helmholtz equation into the following equation, defined in the exterior of the disk $B_{r_0}$:
\begin{align}\label{helmholtzeq2}
\vert (\Phi^{-1})'(x_1+ix_2) \vert^2\partial_t w(t,\bx) - \Delta w + B(v,w) =0,
\end{align}
where is the nonlinear term $(\bv, \nabla w)$ changes to a bilinear form 
$$
B(v,w) = \vert (\Phi^{-1})'(x_1+ix_2) \vert^2(\bv, \nabla w).
$$

Fix a non-negative integer $k$ and divide this equation by $z^k = (x_1+ix_2)^k$. In view
of (\ref{ortcond2}) we get the ratio
\begin{align*}
2\pi r_0^{-|k|}\left(r_0 \frac{\partial w_k(t,r)}{\partial r}\Big|_{r=r_0} + |k| w_k(t,r_0) \right )
+ \int_{\Omega} \frac {B(v,w)}{z^k} d\bx =0,
\end{align*}
which implies a pseudo boundary condition
\begin{align*}
r_0 \frac{\partial w_k(t,r)}{\partial r}\Big|_{r=r_0} + |k| w_k(t,r_0)  = 2 r_0 [G_\tau B(v,w)]_k.
\end{align*}

From this we obtain the energy identity for the domain mapped to the exterior of the disk using the Riemann mapping:
\begin{align*}
\frac 12 \frac d{dt} \left \|w(t,\cdot) \vert (\Phi^{-1})'(z) \vert \right \|^2_{L_2(B_{r_0})} + \nu \|\nabla w\|^2_{L_2(\Omega)} - \nu \|w\|^2_{\dot H^{\frac 12}(\partial \Omega)} + \\ 2 \oint \limits_{|\bx| = r_0}w \cdot G [ B(v,w)] \operatorname{\mathbf{dl}} = 0.
\end{align*}

Consider the integral operator:
\begin{align} \label{BiotSavar}
L[w] = \frac 1{2\pi} D\Phi^t(\bx) \int_{\Omega}   \frac{(\Phi(\bx)-\Phi(\by))^\perp}{\vert \Phi(\bx)-\Phi(\by) \vert ^2} w(\by)d\by,
\end{align}
the kernel of which is the skew gradient of the Green's function:
\begin{align*}
\frac 1{2\pi} \ln |\Phi(\bx)-\Phi(\by)|.
\end{align*}

Here
\begin{align*}
{D\Phi}^t = \begin{pmatrix}
&\frac{\partial y_1}{\partial x_1} &\frac{\partial y_2}{\partial x_1} \\ &-\frac{\partial y_2}{\partial x_1} &\frac{\partial y_1}{\partial x_1}
\end{pmatrix}, 
\end{align*}
and the functions $y_1$, $y_2$ define real valued mapping $\Phi$ in Cartesian coordinates.

Then the energy identity in the $\Omega$ will be written as follows:
\begin{align}\label{energyns}
\frac 12 \frac d{dt} \|w(t,\cdot)\|^2_{L_2(\Omega)} + \nu \|\nabla w\|^2_{L_2(\Omega)} - \nu \|w\|^2_{\dot H^{\frac 12}(\partial \Omega)} + 2 \oint \limits_{\partial \Omega} w \cdot L [ (\bv, \nabla w) ]  \operatorname{\mathbf{dl}} = 0,
\end{align}

In terms of operators $D_k$, defined in (\ref{DkOp}), this identity will take the form:
\begin{align}\label{energy3}
\frac 12 \frac d{dt}  \|w(t,\cdot) \|^2_{L_2(\Omega)} + \nu \sum_{k=-\infty}^\infty \left | D_k[w_k] \right |^2 + 2 \oint \limits_{\partial \Omega}w \cdot L \left [  (\bv, \nabla w) \right ]  \dl = 0.
\end{align}

\medskip
In conclusion, we give an estimate of the last term in the formula (\ref{energy3}). $L_2$-norm of boundary function is controlled by $H^1$-norm. Then with some$C=C(\Omega)>0$ 
$$
\left | \oint \limits_{\partial \Omega}w \cdot L \left [  (\bv, \nabla w) \right ]  \dl \right |
\leq C \|w\|_{H^1(\Omega)} \left \| L \left [  (\bv, \nabla w) \right ] \right \|_{H^1(\Omega)}
$$

This in turn leads to
$$
\left \| L \left [  (\bv, \nabla w) \right ] \right \|_{L_2(\Omega)}
=\left \| \left ( \nabla L , \bv \right ) w  \right \|_{L_2(\Omega)}. 
$$

Operator $\nabla L$ has a singular kernel of the Sigmund-Calderon type and for $p >1$ holds (see \cite{Calderon}\cite{G2}):
\begin{align*}
\left \| \left ( \nabla L , \bv \right ) w  \right \|_{L_p(\Omega)} \leq C \|  w \bv  \|_{L_p(\Omega)}. 
\end{align*} 

The following estimate will also be correct.
\begin{align*}
\left \| \nabla L \left [  (\bv, \nabla w) \right ] \right \|_{L_2(\Omega)} \leq C
\|(\bv, \nabla w)\|_{L_2(\Omega)} \leq C \|\bv\|_{L_\infty(\Omega)}\|\nabla w\|_{L_2(\Omega)} \leq \\
\leq \frac{C \|\bv\|^2_{L_\infty(\Omega)}}{\nu} + \nu \|\nabla w\|^2_{L_2(\Omega)}.
\end{align*}

Then from the energy equality (\ref{energyns}) follows
\begin{align*}
\frac 12 \frac d{dt}  \|w(t,\cdot) \|^2_{L_2(\Omega)} \leq 
\nu \|w\|^2_{\dot H^{\frac 12}(\partial \Omega)} + \frac{C \|\bv\|^2_{L_\infty(\Omega)}}{\nu}.
\end{align*}

This estimate implies the absence of blow-up and the correctness of the external boundary-value problem for 2D Navier-Stokes equations. The bound on $\bv(t,\cdot)$ in $L_\infty(\Omega)$ can be obtained from the inequality:
\begin{align*}
&\|\bv(t,\cdot) - \bv_\infty\|_{L_\infty(\Omega)} \leq C \|\bv(t,\cdot) - \bv_\infty\|^\frac 12 _{L_4(\Omega)} \|\nabla \bv(t,\cdot)\|^\frac 12_{L_4(\Omega)}
\end{align*}

In order to estimate $L_4$-norm of the gradient we can use the estimates of the Bio-Savart operator (see \ref{BiotSavar}). The gradient of the vector field $\nabla \bv(t,\cdot)$ is restored by the vortex function $w$ through an integral transform with a singular kernel of the Sigmund-Calderon type, and the following inequality holds with some $C=C(\Omega)$:
\begin{align*}
&\|\nabla \bv(t,\cdot)\|_{L_4(\Omega)} \leq C \|w(t,\cdot) \|_{L_4(\Omega)}.
\end{align*}

Then considering the inequality
\begin{align*}
&\|w(t,\cdot) \|_{L_4(\Omega)} \leq C \|w(t,\cdot)\|^\frac 12_{L_2(\Omega)} \|\nabla w(t,\cdot)\|^\frac 12_{L_2(\Omega)}
\end{align*}
we will have
\begin{align*}
&\|\nabla \bv(t,\cdot)\|_{L_4(\Omega)} \leq C \|w(t,\cdot)\|^\frac 12_{L_2(\Omega)} \|\nabla w(t,\cdot)\|^\frac 12_{L_2(\Omega)}.
\end{align*}

But $L_4$-norm of the vector field is controlled by $L_2$-norms:
\begin{align*}
\|\bv(t,\cdot) - \bv_\infty\|^4_{L_4(\Omega)} \leq  C \|\bv(t,\cdot)- \bv_\infty\|^2_{L_2(\Omega)} \|\nabla \bv(t,\cdot)\|^2_{L_2(\Omega)}.
\end{align*}

And so
\begin{align*}
&\|\bv(t,\cdot) - \bv_\infty\|_{L_\infty(\Omega)} \leq \\
&~~~~~C_1 \|\bv(t,\cdot)- \bv_\infty\|^\frac14_{L_2(\Omega)} \|\nabla \bv(t,\cdot)\|^\frac14_{L_2(\Omega)} \|w(t,\cdot)\|_{L_2(\Omega)}^\frac 14  \|\nabla w(t,\cdot)\|^\frac 14 _{L_2(\Omega)} \leq \\
&~~~~~~~~~~C_2 \|\bv(t,\cdot)- \bv_\infty\|^\frac14_{L_2(\Omega)} \|w(t,\cdot)\|_{L_2(\Omega)}^\frac 12  \|\nabla w(t,\cdot)\|^\frac 14 _{L_2(\Omega)}.
\end{align*}

\section{Conclusion}

In addition to viscous friction forces, the dynamics of enstrophy is also determined by non-viscous factors. The equation of dynamics of enstrophy (\ref{energyns}) identifies a term that explicitly
 depends on the viscosity  $\nu$, and which slows down the growth of enstrophy. And there is also a term that makes the opposite contribution to the dynamics of enstrophy. This term is determined by the circulation of the vector field
$w \cdot L \left [  (\bv, \nabla w) \right ]$
around the boundary of a streamlined body. At the same time, this value combined with the boundary friction of viscous forces, expressed by the Sobolev norm $\nu \|w\|^2_{\dot H^{\frac 12}(\partial \Omega)}$, distributed in the boundary layer  make a positive contribution to the dynamics. And the term $\nu \|\nabla w\|^2_{L_2(\Omega)}$ distributed over the external domain $\Omega$ provides decrease of enstrophy due to viscosity. For a linearized around zero Stokes system the enstrophy decreases in a power-like manner and
the quadratic form for the vortex function $(\Delta w, w)_{L_2(\Omega)} =  \|\nabla w\|^2_{L_2(\Omega)} - \|w\|^2_{\dot H^{\frac 12}(\partial \Omega)}$ is strictly positive.

\appendix
\section{Numerical simulation}

Here are the results of numerical simulation of vortex flow for the linear Helmholtz equation (\ref{heateq}) and for the nonlinear Helmholtz equation (\ref{helmeq}) with $\nu=1$ and flow velocity $\bv_\infty =(150.0)$. The flow was modeled around a circular and two elliptical cylinders. A no-slip condition was set at the boundary of these bodies.

The flow was calculated in the AGVortex program using the finite element method. The triangulation grid included 67728 elements. The resolution matrix contained 843230 nonzero elements.

\newpage
\begin{figure}[h!]
\center{\includegraphics[scale=0.65]{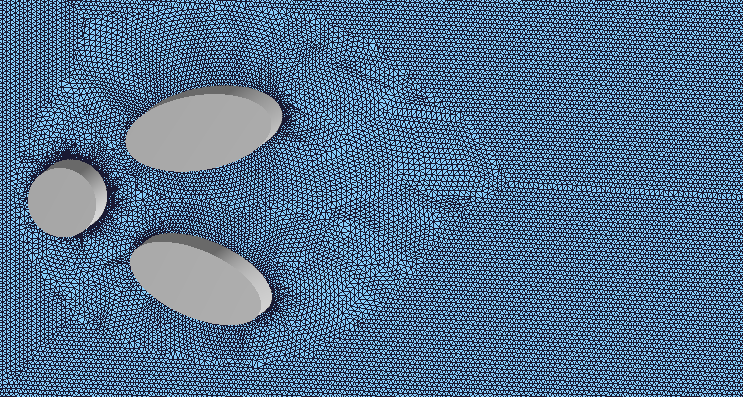}}
\caption{Triangulation grid}
\end{figure}

The exstrophy of the linear Helmholtz equation, which is described by the equation (\ref{energy2}), tends to zero in a power-like manner. The enstrophy of the nonlinear Helmholtz equation is described by (\ref{energy3}). Its dynamics is pseudoperiodic in nature.

\begin{figure}[h!]
\center{\includegraphics[scale=0.65]{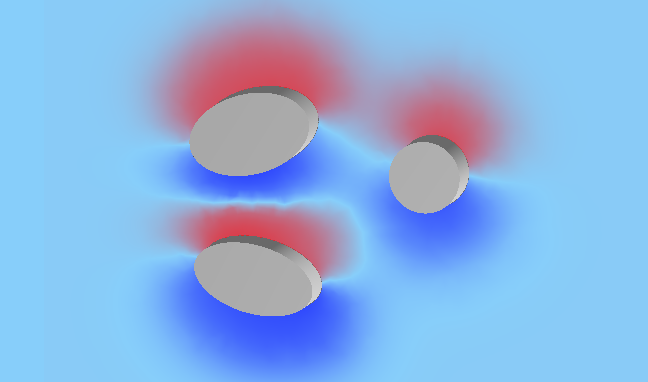}}
\caption{The vortex flow for the Stokes system}
\end{figure}

\begin{figure}[h!]
\center{\includegraphics[scale=0.65]{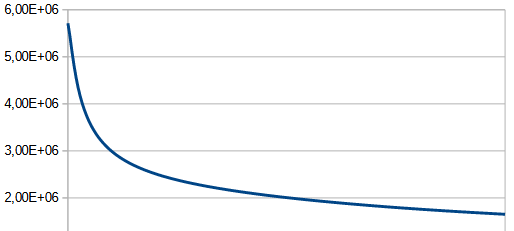}}
\caption{Dissipation of the enstrophy of the Stokes system}
\end{figure}

\begin{figure}[h!]
\center{\includegraphics[scale=0.65, keepaspectratio]{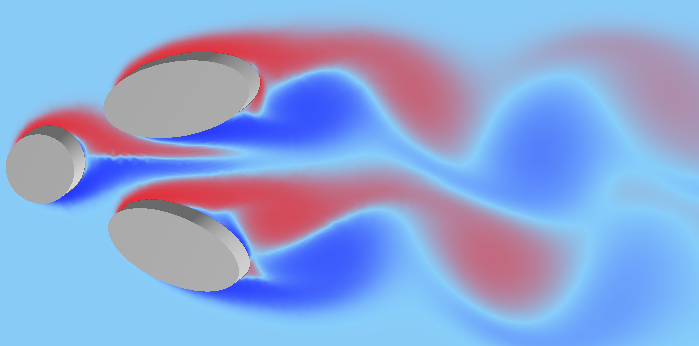}}
\caption{Vortex flow of the Navier-Stokes system around several bodies}
\end{figure}

\begin{figure}[h!]
\center{\includegraphics[scale=0.65]{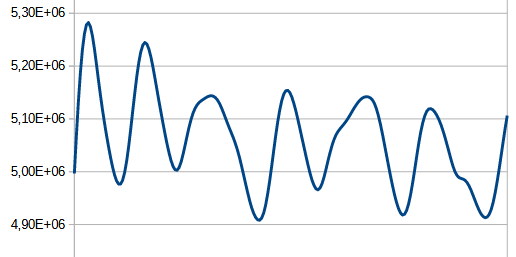}}
\caption{Dynamics of enstrophy for the Navier-Stokes system}
\end{figure}


\newpage


\begin{thebibliography}{9}
\bibitem{BaiesiMaes} Baiesi M, and Maes C. {Enstrophy dissipation in two-dimensional turbulence} // Phys. Rev. E. 2005. V. 72:5. P. 056314. doi:10.1103/PhysRevE.72.056314.

\bibitem{MPY} Matharu P., Protas B., Yoneda T. {On Maximum Enstrophy Dissipation in 2D Navier-Stokes Flows in the Limit of Vanishing Viscosity} // Physica D: Nonlinear Phenomena. 2022. V. 441. P. 133517. doi:10.1016/j.physd.2022.133517.


\bibitem{Wu}Wu J.Z., Wu J.M. {Interactions between a solid surface and a viscous compressible flow field} // J. Fluid Mech. 1993. V. 254 P. 183-211. doi:10.1017/S0022112093002083.

\bibitem{Weiss}Weiss J. {The dynamics of enstrophy transfer in two-dimensional hydrodynamics}// Physica D: Nonlinear Phenomena. 1991. V. 48:2-3. P. 273-294. doi:10.1016/0167-2789(91)90088-Q.

\bibitem{Lighthill} Lighthill J. Introduction: Boundary Layer Theory. Oxford University Press, Oxford, 1963.

\bibitem{Wu2}Wu J.Z. A theory of three-dimensional interfacial vorticity dynamics // Phys. Fluids. 1995. V. 7:10. P. 2375-2395.

\bibitem{Chen} Chen T., Liu T, Wang L.P. Features of surface physical quantities and temporal-spatial evolution
 of wall-normal enstrophy flux in wall-bounded flows // Phys. Fluids. 2021. V. 33:12.  125104.
 
\bibitem{FoiasTemam} Foias, C. Temam, R. Gevrey class regularity for the solutions of the Navier-Stokes
 equations // Journal of Functional Analysis. 1989. V. 87. 359-369.
 
 \bibitem{DoeringGibbon}Doering C. R.,  Gibbon J. D.  Applied Analysis of the Navier-Stokes Equations // 
 Cambridge University Press, New York, 1995.

\bibitem{TL} Tennekes H., Lumley J. L. A First Course in Turbulence. The MIT Press. 1972.


\bibitem{Batchelor}Batchelor G. K. An Introduction to Fluid Dynamics. Cambridge University Press, 2000. 

 
\bibitem{ZhuAntonia}Zhu Y. U, Antonia R.A.  {On the Correlation between Enstrophy and Energy 
Dissipation Rate in a Turbulent Wake} // Appl. Sci. Research. 1997. V. 57. P. 337-347.

\bibitem{Koh}Koh, YM. {Vorticity and viscous dissipation in an incompressible flow} // KSME Journal. 1994. V. 8. P. 35-42. doi:10.1007/BF02953241.

\bibitem{Karman}  Karman, Th. von. {The fundamentals of the statistical theory of turbulence} // J. Aeronaut. Sc. 1937. V. 4:4. P. 131-188. doi:10.2514/8.350. 
 
\bibitem{Taylor}Taylor, G. I. {Production and Dissipation of Vorticity in a 
Turbulent Fluid} // Proc. Roy. Soc. 1938. V. 164:916. P. 15-23.



\bibitem{G2} Gorshkov A. V. {On the unique solvability of the div-curl problem in unbounded domains and energy estimates of solutions} //  Theoret. and Math. Phys., 221:2 (2024), 1799-1812.

\bibitem{G1} Gorshkov A. V. {Special Weber Transform with Nontrivial Kernel} // Mathematical Notes, 114:2 (2023), 172-186.

\bibitem{Lions} Lions, J.L. and Majenes, E. (1971) Inhomogeneous Boundary Value Problems and Their Applications. Mir, Moscow, 371 p.

\bibitem{Calderon} Calderon, A. P., Zygmund, A. {On singular integrals. American Journal of Mathematics}// The Johns Hopkins University Press. 1956. V. 78:2. P. 289-309.

\end{thebibliography}
\end{document}